\newcommand{\be}{\begin{equation}} \newcommand{\ee}{\end{equation}}
\newcommand{\bea}{\begin{eqnarray}}\newcommand{\eea}{\end{eqnarray}}
\newcommand{\nn}{\nonumber}
\newcommand{\ba}{\begin{array}}\newcommand{\ea}{\end{array}}
\newcommand{\R}{\mbox{\bf R}}\newcommand{\U}{\mbox{\bf U}}
\newcommand{\bS}{\mbox{\bf S}}
\newcommand{\tu}{\tilde{u}}\newcommand{\tv}{\tilde{v}}
\newcommand{\tb}{\tilde{b}}
\newcommand\Tq{\mbox{Tr}_q}

\documentstyle[12pt]{article}
\begin{document}
\def\theequation{\arabic{section}.\arabic{equation}}
\renewcommand{\thefootnote}{\fnsymbol{footnote}}
\begin{titlepage}
\rightline {Preprint JINR E2-95-393}
\vspace{2cm}
\begin{center}

 QUANTUM DEFORMATIONS OF THE SELF-DUALITY EQUATION  \\
AND CONFORMAL TWISTORS \\

\vspace{0.5cm}
{\it B.M.Zupnik }\\
{\it Bogoliubov   Laboratory of Theoretical Physics , JINR, Dubna,
Moscow Region, 141980, Russia, E-mail: zupnik@thsun1.jinr.dubna.su }\\

\vspace{0.5cm}
Submitted to Proceedings of VI seminar "Quantum Gravity", Moscow,
June 12-19, 1995
\end{center}

\vspace{1cm}

{\bf Abstract.} A noncommutative algebra of the complex $q$-twistors and
their
differentials is considered on the basis of the quantum $GL_q (4)\times
 SL_q (2)$ group. Real and pseudoreal $q$-twistors are   discussed too.
We consider the quantum-group self-duality equation in the framework
 of the local gauge algebra of differential forms on  $q$-twistor spaces.
 Quantum deformations of the general multi-instanton
solutions are constructed. The corresponding  noncommutative algebras of
moduli  are introduced. The general $q$-instanton connection is a
function of the $q$-twistors and the $q$-moduli .
\end{titlepage}
\renewcommand{\thefootnote}{\arabic{footnote}}
\setcounter{footnote}{0}
\setcounter{section}{0}
\section{Introduction }

$\;\;\;\;$   Noncommutative generalizations of the local gauge theories
have been considered in the framework of different approaches [1-7].
 The mathematically strict approach of Refs[3,4] is
based on the noncommutative global generalizations of the classical
fibre bundles. We prefer to study the local structure of the
quantum-group gauge theory in terms of the deformed connection and
curvature differential forms $A$ and $F$. The basic algebra of these
gauge forms should be covariant under the action of the quantum gauge
group [5,6]
\be
A\;\rightarrow\; T A T^{-1} + dT T^{-1}  \label{A1}
\ee
where $T$ and $dT$ are elements of the differential complex on the
quantum group.

Consider some classical or quantum space with the coordinates $z$ and
let $T(z)$ and $dT(z,dz)$ be noncommutative 'functions' on this space.
We shall treat these functions as generators of the local gauge
differential complex if the map
\be
T\;\rightarrow\; T(z),\;\;\;\;\;dT\;\rightarrow\; dT(z,dz) \label{A2}
\ee
conserves all relations between $T$ and $dT$.

The general function $T(z)$ is some formal expansion with noncommutative
coefficients. Thus, the localization of the quantum group is equivalent
to the definition of an infinite-dimensional noncommutative Hopf algebra.

The $q$-deformations of the Grassmann and twistor spaces were studied in
Refs[8-10]. In  section 2 we consider the
 differential
calculus on the 4-dimensional deformed complex twistor space $T_q(4,C)$.
The real forms of the $q$-twistor space are also discussed .

 Section 3 is devoted to the description of the quantum-group gauge
fields on the $q$-twistor space. We consider the algebraic relation
for the connection form in the $GL_q (N)$ gauge theory [5,6]
that defines the algebraic properties of the 'off-shell' gauge
fields. The generalizations of the Yang-Mills and self-duality
equations are discussed. Note that one can use the 'pure gauge'
$U(1)$ field and $q$-traceless curvature 2-forms in the
$U_q (N)$ gauge theory \cite{a7}. The noncommutative analogue of
the BPST one-instanton solution [11] was constructed in the
 deformed 4-dimensional Euclidean space \cite{a7}.

The deformed analogue of the t'Hooft multi-instanton twistor solution
for the gauge group $GL_q (2)$ (or $U_q (2)$) is considered in
section 4. We use a  multidimensional extension of
the $q$-twistor algebra by the set of  noncommutative $6D$-vector
generators $b$.
The potential of our solution is a sum of the central $(z,b)$-functions
 obeying the $q$-twistor Laplace equation.

The deformed generalization of the Atiyah-Drinfel'd-Hitchin-Manin
solution \cite{a12} for the gauge group $GL_q (N)$ contains $q$-twistor
functions $u$ and $\tu$. We generalize the classical conformal
constructions of Ref\cite{a13}. One can consider the linear twistor
 functions
$v$ and $\tv$ that depend on the noncommutative moduli  $b$ and
$\tb$. The functions $u$ and $v$ are submatrices of the quantum
$GL_q (N+2p)$ matrix for the instanton number $p$. The consistency
 relations for the deformed ADHM-construction can be proved in the
framework of a differential calculus on $GL_q (N+2p)$. The self-duality
condition is equivalent to the bilinear constraint on the moduli
 $b$ and $\tb$.

A preliminary version of this work
was published in Ref[14]. Note that we use here the modified notation and
definitions of some basic quantities.

\setcounter{equation}{0}

\section{Differential calculus on the deformed
\newline twistor space}

$\;\;\;\;$The conformal covariant description of the classical ADHM
solution was considered in Ref\cite{a13}. This approach uses real forms
of the complex $GL(4,C)\times SL(2,C)$ twistors where $GL(4,C)$ is the
complex conformal group. It is convenient to discuss firstly the
deformed complex twistors.

Let $R^{ab}_{cd}$ be a solution of the
Yang-Baxter equation satisfying also the Hecke relation
\bea
& R\;R^\prime \;R = R^\prime \;R\;R^\prime & \\ \label{B1}
& R^2 = I + (q-q^{-1})R  & \label{B2}
\eea
where $q$ is a complex parameter. Note that the standard notation for
these $R$-matrices is $ R=\hat{R}_{12},\;R^\prime=\hat{R}_{23}$
\cite{a15}. We use the symbols $a,b\ldots h=1\ldots 4$ for the 4D-spinor
indices.

The multiparameter 4D $R$-matrix [16-18] and corresponding
inverse matrix $R^{-1}$ can be written
in the following simple form:
\be
(R^{\pm 1})^{ab}_{cd}=\delta^a_c \delta^b_d [q^{\pm 1} - q^{\epsilon (a-b)}]
+ r(ab)\delta^a_d \delta^b_c
\label{B3}
\ee
where $\epsilon (b-a)=0,\;\pm 1$ is a sign function and $r(ab)$ are complex
 parameters
satisfying the relations $r(ab)r(ba)=1$. It is evident that this formula
is valid for an arbitrary number $N$ .

The standard $GL_q (4)$ solution corresponds to the case $r(ab)=1\;$
 \cite{a15}. The choice $q=1$ leads to the unitary $R$-matrix \cite{a19}
\be
R^2 =I \label{B4}
\ee

Consider also the $SL_q(2,C)\;R$-matrix
\be
R^{\alpha\beta}_{\mu\nu}=q\delta^\alpha_\mu \delta^\beta_\nu +
\varepsilon^{\alpha\beta}(q)\varepsilon_{\mu\nu}(q) \label{B5}
\ee
where $\varepsilon(q)$ is the deformed antisymmetric symbol
\be
\varepsilon^{12}(q) =-\varepsilon_{12}(q)=\frac{1}{\sqrt{q}},\;\;\;\;
\varepsilon^{21}(q)=-\varepsilon_{21}(q)=-\sqrt{q}
\label{B6}
\ee

The $q$-deformed flag spaces and twistors were considered in Refs[8-10].
 We shall treat the complex $q$-twistors $z^\alpha_a$ as
generators of the noncommutative algebra with the basic relation
\be
 R^{\alpha\beta}_{\mu\nu}\;z^\mu_a\; z^\nu_b =z^\alpha_c\; z^\beta_d\;
R^{dc}_{ba}  \label{B7}
\ee

This relation for the $(4\times 2)$ rectangular matrix $z$ is analogous
to the RTT-relations for the square quantum matrices. The consistency
conditions for (\ref{B7}) are  pairs of the Yang-Baxter and Hecke
relations (2.1, 2.2) for the independent 4D and 2D $R$-matrices
with the unique common parameter $q$.

The differential calculus on the complex $q$-twistor space $T_q (4,C)$
can be constructed by the analogy with the bicovariant differential
complex on the quantum linear group [20-23]. Consider
the relations between $z^\alpha_a$ and their differentials $dz^\alpha_a$
\bea
& z^\alpha_a\;dz^\beta_b = R^{\alpha\beta}_{\mu\nu}\;dz^\mu_c\; z^\nu_d\;
R^{dc}_{ba} &\\ \label{B8}
& dz^\alpha_a\;dz^\beta_b = -R^{\alpha\beta}_{\mu\nu}\;dz^\mu_c\; dz^\nu_d
\;R^{dc}_{ba} & \label{B9}
\eea

The elements $z$ and $dz$ are generators of the external algebra
$ \Lambda T_q (4,C)$ on the $q$-twistor space. The operator of external
derivative $d$ on $ \Lambda T_q (4,C)$ is nilpotent and satisfies the
ordinary Leibniz rule.

The symmetry properties of $dz$ can be obtained from Eq(\ref{B9})
\be
P^{(+)}_2 dz\;dz^\prime P^{(+)}_4 = 0 = P^{(-)}_2 dz\;dz^\prime P^{(-)}_4
\label{B10}
\ee
where $P^{(\pm)}_2$ and $ P^{(\pm)}_4$ are the projection operators for
$SL_q (2)$ and $GL_q (4)$,  respectively \cite{a15}
\be
P^{(+)} + P^{(-)}=I,\;\;\;\;R =qP^{(+)} - q^{-1}P^{(-)}
\label{B11}
\ee

One can define the algebra of partial derivatives  $\partial^a_\alpha$
on $T_q (4,C)$
\bea
& R^{ab}_{cd}\;\partial^c_\alpha\; \partial^d_\beta=\partial^a_\mu\;
\partial^b_\nu\; R_{\beta\alpha}^{\nu\mu} & \\ \label{B12}
& \partial^a_\alpha\; z^\beta_b =\delta^a_b\; \delta_\alpha^\beta +
R^{\beta\mu}_{\alpha\nu}\; R^{da}_{cb}\; z^\nu_d\; \partial^c_\mu &\\
 \label{B13}
&  \partial^a_\alpha\; dz^\beta_b = R^{\beta\mu}_{\alpha\nu}\;
R^{da}_{cb}\; dz^\nu_d\; \partial^c_\mu
\label{B14}
\eea

Consider a definition of the deformed $\varepsilon_q$-symbol
for $GL_q (4)$
\be
\varepsilon_q^{abcd}=-q R^{ba}_{fe}\; \varepsilon_q^{efcd}=
[P_4^{(-)}]^{ba}_{fe}\; \varepsilon_q^{efcd}=-q^{\epsilon (b-a)}r(ba)\;
\varepsilon_q^{bacd}
                     \label{B15}
\ee
 Analogous relations are valid for other neighboring pairs of indices.

The $q$-twistors obey the following identity:
\be
\varepsilon_q^{abcd}\; z^\beta_b\; z^\mu_c\; z^\nu_d = 0 \label{B16}
\ee

Introduce the  $SL_q(2)$-invariant bilinear function of $q$-twistors
\be
 y_{ab}=\frac{q^2}{1+q^2}\;\varepsilon_{\alpha\beta}(q)\; z^\alpha_a\;
 z^\beta_b =[P^{(-)}_4]^{dc}_{ba}\;y_{cd}  \label{B17}
\ee
This vector satisfies the following commutation relation
\be
y_{ab}\;z^\alpha_c =q^{-1} R^{ed}_{ha} R^{fh}_{cb}\;z^\alpha_d\; y_{ef}
\label{B18}
\ee

In consequence of Eq(\ref{B16}) the coordinate $y$ is an isotropic
vector in the deformed 6D space
\be
 (y,y)=\varepsilon_q^{abcd}\; y_{ab}\; y_{cd}=0  \label{B19}
\ee
This $GL_q (4,C)$ covariant equation determines the 4D deformed subspace
of the complex 6D quantum plane. The classical analogue of this subspace
is a complex sphere $S_4^C$.

Consider a duality transformation $\ast$ of the basic $q$-twistor 2-forms
by analogy with Ref\cite{a7}
\be
\ast dz\; dz^\prime = dz\; dz^\prime\; P^{(+)}_4 - dz\; dz^\prime\;
P^{(-)}_4
\label{B20}
\ee
where $P^{(\pm)}_4$ are the  $GL_q(4)$ projectors (\ref{B11}).
Note that a self-dual part $dz\; dz^\prime\; P^{(+)}_4$ is proportional
 to the $SL_q (2)$-invariant conformal tensor
\be
\varepsilon_{\alpha\beta}(q)\; dz^\alpha_a\; dz^\beta_b \label{B21}
\ee

Real $q$-twistors can be treated as a representation of the real quantum
 group
$SL_q (2,R)\times GL_q (4,R)$. The classical analogue of these twistors
are connected with the real pseudo-Euclidean $(2,2)$-space \cite{a13}.
Consider the $R$-matrices (\ref{B3},\ref{B5}) and the conditions
$|q|=1$ and $|r(ab)|=1$, then under the complex conjugation
\bea
& \overline{R^{ab}_{cd}}=(R^{-1})^{ba}_{dc} &\\ \label{B22}
& \overline{R^{\alpha\beta}_{\mu\nu}}=(R^{-1})^{\beta\alpha}_{\nu\mu}
 & \label{B23}
\eea

These formulas correspond to anti-involution of
 the real  $T_q (2,2)$ twistors
\bea
& \overline{z}=z,\;\;\;\;\;\overline{dz}=dz &\\ \label{B24}
& \overline{z\;z^\prime}= z^\prime\;z,\;\;\;\; \overline{z\;dz^\prime}=
dz^\prime\;z & \label{B25}
\eea

The pseudoreal Euclidean $q$-twistors have a more complicated
anti-involution
\be
\overline{z^\alpha_a}=\varepsilon_{\alpha\beta}(q) z^\beta_b\; C^b_a (q)
\label{B26}
\ee
where $C(q)$ is the charge-conjugation matrix for the Euclidean
conformal quantum group $U^*_q (4)=  D \times SU^*_q (4)$ and $D$ is a
 real one-parameter dilatation.

It is convenient to use the simple representation
\be
C(q)=\left(\begin{array}{cc} \varepsilon^{\alpha\beta}(q)& 0 \\
		 0 & \varepsilon^{\dot{\alpha}\dot{\beta}}(q)\\
		\end{array} \right)  \label{B27}
\ee

For the real $q$ we have the following properties:
\bea
& \overline{\varepsilon^{\alpha\beta}(q)}=\varepsilon^{\alpha\beta}(q)=
-\varepsilon_{\alpha\beta}(q) &\\ \label{B28}
& \overline{C(q)}=C(q),\;\;\;\;\;C^2 (q)=-I &\\ \label{B29}
& \overline{\overline{z^\alpha_a}}=-\varepsilon^{\alpha\beta}(q)
\varepsilon_{\beta\gamma}(q) z^\gamma_b\; (C^2 )^b_a (q)=z^\alpha_a &\\
\label{B30}
&\overline{R^{dc}_{ba}} =C^d_e (q)\; C^c_f (q)\; R^{fe}_{gh}\; C^g_a (q)
 \;C^h_b (q) & \label{B31}
\eea

We do not here consider the conformal quantum group $SU_q (2,2)$ [24,25]
and corresponding $q$-twistors.

\setcounter{equation}{0}

\section{Quantum-group gauge theory on the \newline $q$-twistor space}

$\;\;\;\;$The classical gauge field on some domain $\{x^m \}$ of the
basic space corresponds to the connection 1-form
\be
A(x,dx)=dx^m A_m (x) \label{C1}
\ee
which can be decomposed in terms of the gauge-group generators. For the
domain with the coordinates $\tilde{x}(x)$ one should define the
transformed connection
\be
\tilde{A}(\tilde{x},d\tilde{x})= T(x) A T^{-1}(x) + dT(x) T^{-1}(x)
\label{C2}
\ee
where $ T(x)$ is a matrix of the local gauge transformation.

The components of the matrix $A(x,dx)$ satisfy the anticommutativity
conditions
\be
\{A^i_k , A^l_m \} =0 \label{C3}
\ee

The classical gauge group formally has an infinite number of generators.
A constructive example of the classical gauge algebra is the affine
(Kac-Moody) algebra. The quantum affine algebras can be considered as a
basis of the quantum gauge theory on the classical two-dimensional space.

The formal quantization of the gauge groups on the multi-dimensional
classical or quantum spaces is a difficult problem.  Let $R_N $ be the
constant $R$-matrix for the quantum group $GL_q (N)$ and $x^M$ are the
coordinates of some basic space. Consider the simplest possible relations
for the components of the quantum gauge matrix
\bea
& R_N T(x) T^\prime (x)=T(x) T^\prime (x) R_N &\\ \label{C4}
& T^i_k (x)=\sum\limits^\infty_{n=0} \frac{1}{n!} (T^i_k )_{M_1\cdots M_n}
x^{M_1}\ldots x^{M_n} & \label{C5}
\eea
where $i,k\ldots =1\ldots N$.

The quantum-group gauge matrix is well defined if the relation (3.4)
generates the consistent set of relations between the coefficients
$(T^i_k )_{M_1\cdots M_n} $.

Quantum deformations of the $GL(N)$ gauge connection can be treated in
 terms of the noncommutative gauge algebra for the components of the
 deformed connection 1-form [5,6]
$$
(A R_N A + R_N A R_N A R_N )^{ik}_{mn} =
$$
\be
 A^i_l (R_N )^{lk}_{rn} A^r_m +
(R_N )^{ik}_{jl} A^j_r (R_N )^{rl}_{st} A^t_p (R_N )^{pt}_{mn} = 0
\label{C6}
\ee
These relations generalize the classical anticommutativity conditions
(\ref{C3}). The gauge algebra is an analogue of the relations between
components of the right-invariant 1-forms $\omega =dT T^{-1} $ in the
framework of the bicovariant differential calculus on $GL_q (N)$ [20-23]
\be
\omega R_N \omega + R_N \omega R_N \omega R_N = 0,\;\;\;d\omega -
\omega^2 = 0 \label{C7}
\ee

Thus, the form $\omega$ can be considered as a pure gauge $GL_q (N)$-field.
The general $GL_q (N)$ connection $A$ has the nontrivial curvature 2-form
\be
F=dA - A^2 \label{C8}
\ee

 Explicit constructions of the deformed gauge fields on the $q$-twistor
space contain also the noncommutative elements ( moduli) which generate
 some algebra $B$
\be
A^i_k (z,dz,B) = dz^\alpha_a (A^a_\alpha )^i_k (z,B) \label{C9}
\ee
The appearance of  additional noncommutative elements is necessary for
the consistency of the algebra (\ref{C6}) and the relations of the
$q$-twistor algebra (2.8,2.9).

The (anti)self-duality equation for the gauge field (\ref{C9}) can be
defined with the help of the relations (\ref{B20})
\be
\ast F =\frac{1}{2}(\ast dz^\alpha_a dz^\beta_b ) F^{ab}_{\alpha\beta}(z,B)
=\pm F \label{C10}
\ee
or in terms of the deformed field-strength
\be
[P^{(\pm)}_4 ]^{ab}_{cd} F^{cd}_{\alpha\beta}(z,B) = 0 \label{C11}
\ee

The solutions of this equation and the explicit construction of the
 algebra $B$ will be considered in  sections 4 and 5.

A quantum deformation of the Yang-Mills equation has the standard form in
the framework of the external algebra  $\Lambda T_q (4,C)$
\be
\nabla \ast F = d \ast F + [A,\ast F ] = 0 \label{C12}
\ee

The bicovariant differential calculus with the ordinary Leibniz rule for
the operator $d$ and the gauge-connection algebra (\ref{C6}) are consistent
only for the case of the nonsemisimple quantum group $GL_q (N)$. The gauge
algebra produce the restriction
\be
\alpha=\Tq A \ne 0 \label{C13}
\ee

Nevertheless, one should use the gauge-covariant conditions \cite{a7}
\be
d\alpha =0,\;\;\;\Tq A^2 =0,\;\;\;\Tq F =0 \label{C14}
\ee
These restrictions generate the effective reduction of the Abelian gauge
field $\alpha$ in the framework of the gauge group $GL_q (N)$.

\setcounter{equation}{0}

\section{Quantum deformations of the t'Hooft \newline multi-instanton
 solution}

$\;\;\;\;$ A simple form of the manifest  multi-instanton solution in
the Euclidean space was discussed in Refs[26-28]. This solution can be
written in terms of the potential $\Phi $ satisfying the Laplace equation.
The classical twistor version of the t'Hooft solution was considered in
Ref[13].

Consider firstly the deformed Laplace equation in the complex $q$-twistor
space $T_q (4,C)$. Using Eq(2.13) one can obtain the action of the
 partial\\ $q$-twistor derivative on the isotropic $6D$-vector
\be
\partial^c_\alpha\; y_{ab} =\varepsilon_{\alpha\beta}(q)
 \;[P^{(-)}_4]^{dc}_{ba}\; z^\beta_d
\label{D1}
\ee

Introduce the formal differential operator that  acts only on the $6D$
 vector variables
\be
\partial^{dc}\triangleright y_{ab} = [P^{(-)}_4]^{dc}_{ba}
\label{D2}
\ee

Now we can write the following relations:
\bea
& dy_{ab}=\varepsilon_{\alpha\beta}(q)\;[P^{(-)}_4]^{dc}_{ba}\;dz^\alpha_c
\;z^\beta_d &\\ \label{D3}
& d\Phi (y)=dz^\alpha_c\;\partial^c_\alpha\; \Phi (y)=dy_{ab}\;
\partial^{ba}\;\Phi (y)& \\ \label{D4}
& \partial^c_\alpha\;\Phi (y) =\varepsilon_{\alpha\beta}(q)\;z^\beta_b\;
\partial^{bc}\; \Phi (y)& \label{D5}
\eea

The $SL_q (2,C) $-invariant analogue of the Laplace operator has the
following form:
\be
\Delta^{ba} =-\frac{q}{1+q^2}\varepsilon^{\alpha\beta}(q)\;\partial^b_\beta
\;\partial^a_\alpha                \label{D6}
\ee
By definition, we have
\be
\Delta^{ba} \triangleright y_{cd}=[P^{(-)}_4]_{dc}^{ba}=\partial^{ba}
 \triangleright y_{cd} \label{D7}
\ee

In this section we shall use the standard $GL_q (4,C)\;\;R$-matrix
corresponding to Eq(\ref{B3}) with $r(ab)=1$. This $R$-matrix satisfies
the following identity:
\be
\varepsilon_q^{abcd}\; R^{a^\prime h}_{ea}\; R^{b^\prime e}_{fb}\;
R^{c^\prime f}_{gc}\; R^{d^\prime g}_{h^\prime d}=q\;\delta^h_{h^\prime }
\;\varepsilon_q^{a^\prime b^\prime c^\prime d^\prime}
\label{D8}
\ee

Introduce the additional noncommutative moduli  $b^p_{ab}$
where $p$ is an arbitrary  number
\bea
&(b^p , b^p )= \varepsilon_q^{abcd} b^p_{ab}\; b^p_{cd} = 0 &\\ \label{D9}
& y_{ab}\; b^p_{cd} = R^{ea^\prime}_{ga}\; R^{fg}_{cb}\; R^{c^\prime
b^\prime}_{he}\; R^{d^\prime h}_{df}\; b^p_{a^\prime b^\prime}\;y_{c^\prime
d^\prime} &\\ \label{D10}
& b^p_{ab}\; b^{\hat{p}}_{cd} =q^{-2} R^{ea^\prime}_{ga}\; R^{fg}_{cb}\;
R^{c^\prime b^\prime}_{he}\; R^{d^\prime h}_{df}\;b^{\hat{p}}_{a^\prime
 b^\prime}\;b^p_{c^\prime d^\prime} & \label{D11}
\eea
and $p \le {\hat{p}}$ in the last equation.

This $(B,y)$-algebra has the following central elements
\be
X_p = (y, b^p )= \varepsilon_q^{abcd}\; y_{ab}\; b^p_{cd} \label{D12}
\ee
The commutativity of $X_p$ with $y$ and $b^{\hat{p}}$ can be proved with
 the help of Eq(\ref{D8}).

Let us introduce the commutation relations between $b^p,$ and $z$
\be
b^p_{ab}\; z^\gamma_c = R^{eh}_{ga}\; R^{fg}_{cb}\; z^\gamma_h\;
b^p_{ef}  \label{D13}
\ee
 An analogous relation for
$b^p,$ and $dz$ can be obtained as the external derivative $d$ of this
 formula by using $ [d, b^p ] = 0$.

Equation(\ref{B8}) generates the relation for $y$ and $dz$
\be
y_{ab}\;dz^\gamma_c = q\; R^{eh}_{ga}\; R^{fg}_{cb}\;dz^\gamma_h\; y_{ef}
\label{D14}
\ee

Write the corresponding relation for the elements $X_p$
\be
X_p\; dz^\alpha_a = q^2\; dz^\alpha_a\; X_p \label{D15}
\ee

Now one can determine the derivative of the central functions
\be
\partial^a_\alpha\; \frac{1}{X_p} =-   \frac{ 1 }{ q^2 X_p^2 }  \;
 \partial^a_\alpha X_p  \label{D16}
\ee
\be
\partial^a_\alpha\; \frac{1}{X_p^2}   =-\frac{1+q^2}{q^4}\frac{1}{X_p^3}
\partial^a_\alpha X_p  \label{D17}
\ee

It is easy to check the following identity for the  isotropic
 vectors $b^p$ :
\be
\varepsilon^{\alpha\beta}(q)\;\partial^b_\beta X_p \partial^a_\alpha X_p =
-q^{-1}\varepsilon_q^{abcd}\; b^p_{cd}\; X_p \label{D18}
\ee

We can obtain the solutions of the deformed Laplace equation\\ ($q$-harmonic
functions )
\be
\Delta^{ba}\;\frac{1}{X_p} = \frac{1}{q^5 (1+q^2 )}\frac{1}{X_p^2}
\varepsilon^{\alpha\beta}(q)\;[\partial^b_\beta \;\partial^a_\alpha X_p -
(1+q^2)\frac{1}{X_p}\partial^b_\beta X_p \partial^a_\alpha X_p ]=0
\label{D19}
\ee

By analogy with Ref[13] one can consider the deformed t'Hooft Ansatz
for the $GL_q (2)$ self-dual gauge field
\bea
&A^\alpha_\beta =q^{-3}\;dz^\alpha_a\; (\partial^a_\mu \Phi)\;\Phi^{-1}
\varepsilon^{\sigma\mu}(q)\varepsilon_{\sigma\beta}(q)&\\ \label{D20}
& \Tq A =-q^3 d\Phi\; \Phi^{-1},\;\;\;\; \Tq dA = 0 & \label{D21}
\eea
where the potential function $\Phi$ for the instanton number $P$ is a
 sum of $q$-harmonic functions
with  different elements $b^p$
\be
\Phi = \sum_{p=1}^P \frac{1}{X_p} =\sum_{p=1}^P (y,b^p )^{-1} \label{D22}
\ee

The anti-self-dual part of the corresponding curvature form vanishes in
consequence of Eq(\ref{D19})
\be
(F -\ast F)^\alpha_\beta\; \sim\; dz^\alpha_e\; dz^\gamma_f\;
 [P^{(-)}_4 ]^{fe}_{ac}\; \Delta^{ac} \Phi\;
\Phi^{-1}\;\varepsilon_{\beta\gamma}(q) = 0  \label{D23}
 \ee

Note that the isotropic vector $b^p$ has 5 independent elements so
(4.20) is the $5P$-parameter solution.

\setcounter{equation}{0}

\section{ Quantum deformations of the ADHM-solution }

$\;\;\;\;$The covariant formulation of the ADHM multi-instanton solution
in the classical twistor space was considered in Ref[13]. We shall discuss
the quantum deformations of this formalism.

Let us consider the gauge group $GL_q (N,C)$. The ADHM-solution for the
instanton number $p$ can be connected with some $GL_q (N+2p,C)$ matrix
$q$-twistor function. Introduce the notation for indices of different
types : $A,B\ldots=1\ldots p;\;\;I,K,L,M\ldots =1\ldots N+2p $ and
$i,k,l\ldots =1\ldots N$. The ADHM Ansatz for the general self-dual
$GL_q (N,C)$ field contains the deformed twistor functions $u^i_I (z)$
and $\tu^I_i (z)$
\be
A^i_k =du^i_I\; \tilde{u}^I_k,\;\;\;\; u^i_I\; \tilde{u}^I_k=\delta^i_k
\label{E1}
\ee

The commutation relations for the $u$ and $ \tilde{u}$ twistors are
\bea
& (R_N)^{ik}_{lm}\;u^l_I\; u^m_K = u^i_L\; u^k_M\; \R^{LM}_{IK} &\\
\label{E2}
& \R^{KI}_{ML}\;\tilde{u}^L_i\; \tilde{u}^M_k =\tilde{u}^I_l\;\tilde{u}^K_m
 \;(R_N)_{ki}^{ml} & \\ \label{E3}
& \tilde{u}^I_l\; (R_N)^{li}_{mk}\; u^m_K = u^i_L\;
\R^{IL}_{KM}\;\tilde{u}^M_k & \label{E4}
 \eea
 where the $R$-matrices for $GL_q(N,C)$ and $GL_q(N+2p,C)$ are used.

Introduce also the relation for the differentials $du$
\bea
& \tilde{u}^I_i\;(R_N )^{ik}_{lm}\;du^l_K = du^k_L\;(\R^{-1})^{IL}_{KM}\;
\tilde{u}^M_m &\\ \label{E5}
& du^i_L\;du^k_M\;(\R^{-1})^{LM}_{IK}= -(R^{-1}_N )^{ik}_{lm}\;
du^l_I\;du^m_K & \label{E6}
\eea

These relations are necessary for proving a validity of the gauge algebra
 (\ref{C6}) in the framework of the  ADHM-Ansatz
\be
(A R_N A)^{ik}_{mn}=du^i_I\; du^k_L\; (\R^{-1})^{IL}_{KM}\;\tu^M_n\;
 \tu^K_m = -(R_N A R_N A R_N )^{ik}_{mn} \label{E7}
\ee

Consider also the linear twistor functions $v$ and $\tv$
\bea
& v^{A\alpha}_I = z^\alpha_a \;b^{aA}_I &\\ \label{E8}
& \tilde{v}^{IA\alpha} =z^\alpha_a\; \tilde{b}^{aIA}  & \label{E9}
\eea
where $b$ and $\tb$ are the noncommutative $q$-instanton  moduli
\bea
&b^{aA}_I\; z^\alpha_b = R^{da}_{cb}\; z^\alpha_d\; b^{cA}_I &\\\label{E10}
&\tilde{b}^{aIA}\; z^\alpha_b = R^{da}_{cb}\; z^\alpha_d\; \tb^{cIA} &
 \label{E11}
\eea
The relations between $b$ and $\tb$ will be defined below.

Introduce the following condition for the functions  $v$ and $\tv$:
\be
v^{A\alpha}_I\; \tilde{v}^{IB\beta}=g^{AB}(z)\;\varepsilon^{\alpha\beta}(q)
\label{E12}
\ee
where $g(z)$ is the nondegenerate $(p\times p)$ matrix with the central
elements
\be
g^{AB}(z)=q^{-2}\; y_{cd}\; b^{cA}_I\; \tilde{b}^{dIB}  \label{E13}
\ee

The condition (\ref{E12}) is equivalent to the restriction on the
elements of the $B$-algebra
\be
[P^{(+)}]^{ab}_{cd}\; b^{cA}_I\; \tilde{b}^{dIB} =0 \label{E15}
\ee

Write the basic commutation relations of the $B$-algebra
\bea
& R^{ab}_{cd}\; b^{cA}_I\; b^{dB}_K =b^{aB}_L\; b^{bA}_M\; \R^{ML}_{KI}&\\
\label{E16}
& \R_{LM}^{IK}\;\tilde{b}^{aLA}\;\tilde{b}^{bMB}=R^{ab}_{cd}\;
\tilde{b}^{cIB}\;\tilde{b}^{dKA} &\\ \label{E17}
& R^{ab}_{cd}\;b^{cA}_I\; \tilde{b}^{dKB}=\R^{KL}_{IM}\;\tilde{b}^{aMB}\;
b^{bB}_L
& \label{E18}
\eea
Remark that a formal permutation of the indices $A$ and $B$ is
commutative.

Consider the new functions
\be
\tilde{v}_{A\alpha}^I=
\tilde{v}^{IB\beta}\;g_{BA}(z)\;\varepsilon_{\beta\alpha}(q)  \label{E19}
\ee
where we use the  matrix $g_{BA}$ inverse of the matrix (\ref{E13})
\be
g_{BA}(z)\; g^{AC}(z) = \delta^C_B \label{E20}
\ee

Now one can construct the full quantum $GL_q(N+2p,C)$ matrices
\be
\U=\left( u^i_I\atop v^{A\alpha}_I \right)\;,\;\;\;\;\;\;
\bS(\U)=\U^{-1}=\left(
\tilde{u}^I_i \atop\tilde{v}_{A\alpha}^I \right)
\label{E21}
\ee

The standard $GL_q(N+2p,C)$ commutation relations for these matrices
contain Eqs(5.2-5.4) and the relations for $v$ and $\tilde{v}$
 functions
\bea
& \widetilde{\R}\; \U\; \U^\prime = \U\; \U^\prime\; \R &\\ \label{E22}
& \R\; \bS^\prime\; \bS =\bS^\prime\; \bS\; \widetilde{\R} &\\ \label{E23}
& \bS\; \widetilde{\R}\; \U = \U^\prime\; \R\; \bS^\prime & \label{E24}
\eea
where the $R$-matrix for $GL_q (N+2p,C)$ can be written in the following
form
$$
\widetilde{\R}= \left( \begin{array}{cccc}
(R_N )^{ik}_{mn} &   0    &    0    &  0  \\
 0  & \lambda\delta^A_C \delta^\alpha_\mu \delta^k_n &  \delta^A_D
\delta^\alpha_\nu \delta^k_m &  0 \\
 0 & \delta^i_n \delta^B_C \delta^\beta_\mu & 0  & 0 \\
 0 & 0  & 0 & \delta^A_D \delta^B_C R^{\alpha\beta}_{\mu\nu} \end{array}
\right)
$$
where $\lambda =q - q^{-1} $.

The equations (5.5,5.6) follow from are the relations
\be
\bS\; \widetilde{\R}\; d\U = d\U^\prime\; \R^{-1}\; \bS^\prime \label{E25}
\ee

It should be stressed that the bicovariant differential calculus on
$GL_q (N+2p,C)$ is the basis of the deformed ADHM-construction for the
group $GL_q (N,C)$.

Write explicitly the orthogonality and completeness conditions for the
deformed ADHM-twistors:
\bea
& u^i_I\;\tilde{v}^{IA\alpha} =0 &\\ \label{E26}
& v^{A\alpha}_I \;\tilde{u}^I_i=0 &\\ \label{E27}
& \delta^I_K =\tilde{u}^I_i\;u^i_K\;+\;\tilde{v}^{IA\alpha}\;
g_{AB}(z)\varepsilon_{\alpha\beta}(q)\;v^{B\beta}_K &\label{E28}
\eea

Now we are in a position to verify the self-duality of the connection
(\ref{E1})
\bea
& dA^i_k - A^i_l\;A^l_k = du^i_I\;(\tilde{u}^I_l\;u^l_M\;-\;\delta^I_M)
d\tilde{u}^M_k = &\\ \nn
& =-q^{-4}u^i_I\;\tilde{b}^{cIA} D^a_c\;g_{AB}(z)\varepsilon_
{\alpha\beta}(q)\;
dz^\alpha_a\;dz^\beta_b\;b^{bB}_M\;\tilde{u}^M_k & \label{E29}
\eea
where $D^a_c$ is the $GL_q (4)$ metric.
This curvature contains  only the self-dual $q$-twistor 2-form
(\ref{B21}).

The real forms of the deformed ADHM-construction are based on the quantum
groups $U_q (N)$ and $GL_q (N,R)$.

It should be stressed that all $R$-matrices of our deformation scheme
satisfy the Hecke relation with the common parameter $q$. The other
possible parameters of different $R$-matrices are independent. The case
$q=1$ corresponds to the unitary deformations $(R^2=I) $ of the twistor
space and the gauge groups. It is evident that the trivial deformation
of the $z$-twistors is consistent with the nontrivial unitary deformation
of the gauge sector and vice versa.

The author would like to thank
A.T. Filippov, E.A. Ivanov, A.P.Isaev and
V.I. Ogievetsky
 for helpful discussions and interest in this work.

I am grateful to the administration of JINR and Laboratory of Theoretical
Physics for hospitality. This work was supported in part by
the ISF-grant RUA000, INTAS-grant 93-127  and the contract No.40 of Uzbek
 Foundation of Fundamental
Research .

\end{document}